\documentclass[11pt]{article}

\usepackage{amssymb}
\usepackage{amsthm}
\usepackage{eucal}
\usepackage{amsmath,amsxtra}

\makeatletter

\skip\footins 24pt plus 4pt minus 2pt

\makeatother

\newcommand\A{\mathcal A}
\newcommand\D{\mathcal D}
\newcommand\E{\mathcal E}
\newcommand\Hs{\mathcal H}

\renewcommand\L{\mathcal L}

\newcommand\s{\mathcal S}

\newcommand{\CC}{\mathbf{C}}

\newcommand{\GG}{\mathbf{G}}
\newcommand{\HH}{\mathbf{H}}
\newcommand{\MM}{\mathbf{M}}
\newcommand{\OO}{\mathbf{O}}
\newcommand{\RR}{\mathbf{R}}
\newcommand{\SL}{\mathbf{S}\mathbf{L}}
\newcommand{\SU}{\mathbf{S}\mathbf{U}}
\renewcommand{\SS}{\mathbf{S}}
\newcommand{\ZZ}{\mathbf{Z}}

\newcommand\Cinf{\mathcal{C}^\infty}
\renewcommand{\Im}{\mathop{\rm Im\,}\nolimits}

\newcommand{\thh}{\mathop{\rm th}\nolimits}

\newcommand\bigcheck[1]{#1 \raise1ex\hbox{$\hspace{-1ex}{}^\vee$}}
\newcommand\sucheck[1]{#1 \raise0.5ex\hbox{$\hspace{-1ex}{}^\vee$}}

\newcommand\subH{{\raise-1ex\hbox{$\scriptstyle H$}}}

\newcommand\barsubH{{\raise-1ex\hbox{$\left| \scriptstyle H \right|$}}}
\newcommand\barsubHo{{\raise-1ex\hbox{$\left| \scriptstyle H_0 \right|$}}}

\makeatletter
\@addtoreset{equation}{section}
\makeatother

\newtheorem{theorem}{Theorem}[section]

\renewenvironment{proof}[1][Proof]
           {\medbreak\noindent \emph{#1: \enspace}}

\newenvironment{remark}[1][Remark]
           {\medbreak\noindent \textbf{#1 \enspace}}
           {\par \medbreak}

\makeatletter
\def\iint{\DOTSI\protect\ints@\tw@}
\def\iiint{\DOTSI\protect\ints@\thr@@}
\def\iiiint{\DOTSI\protect\ints@{4}}
\def\idotsint{\DOTSI\protect\ints@\z@}

\def\intkern@{\mkern-6mu\mathchoice{\mkern-3mu}{}{}{}}
\let\DOTSI\relax
\let\ilimits@\displaylimits

\def\ints@#1{%
  \mkern-7mu\mathchoice{\mkern-2mu}{}{}{}%
  \mathop{\mkern7mu\mathchoice{\mkern2mu}{}{}{}%
    \intop\ifnum#1=\z@\intdots@
    \else\intkern@\fi
    \ifnum#1>\tw@\intop\intkern@\fi
    \ifnum#1>\thr@@\intop\intkern@\fi
    \intop
  }\ilimits@
}
\makeatother

\newcommand{\romanparenlist}{
  \renewcommand{\theenumi}{\roman{enumi}}%
  \renewcommand{\labelenumi}{(\theenumi)}%
}

\newcommand{\st}[1]{\ensuremath{^{\scriptstyle \textrm{#1}}}}

\makeatletter
\def\@maketitle{\newpage
 \null
 \vskip 2em
 \begin{center}%
  {\Large\bf \@title \par}%
  \vskip 1.5em
  {\normalsize
   \lineskip .5em
   \begin{tabular}[t]{c}\@author
   \end{tabular}\par}%
  \vskip 2em
  {\@date}%
 \end{center}%
 \par
 \vskip 2.5em}
\makeatother

\begin{document}

\hfill J\'anos Bolyai Memorial Conference\\
\hspace*{3.35in} Budapest, July 2002\\

\vspace{4ex}

\begin{center}
\Large{\textbf{Non-Euclidean Analysis}}
\end{center}
\begin{center}
{\sc SIGURDUR HELGASON}
\end{center}

\begin{center}
June 30, 2003
\end{center}



\section{The Non-Euclidean Plane}
\label{sec:1}

In case the work of Bolyai \cite{Bo} and Lobatschevsky \cite{Lo} left
any doubts about the existence of non-Euclidean geometry these
doubts were removed by the work \cite{Be} of Beltrami.  With a
modification made  possible by hindsight one can state the
following result.

\begin{theorem}
  \label{th:1.1}

Given a simply connected region $D \subset \RR^2 \, (D \neq \RR^2)$ there exists a
Riemannian metric $g$ on $D$ which is invariant under all
conformal transformations of $D$.  Also $g$ is unique up to a
constant factor.
\end{theorem}

Because of the Riemann mapping theorem we can assume $D$ to be the
unit disk.    Given $a \in D$ the mapping $\displaystyle{\varphi
  :z \to \frac{a-z}{1-\bar{a}z}}$ is conformal and $\varphi (a)
=0$.  The invariance of $g$ requires
\begin{equation}
  \label{eq:1.1}
  g_a (u,u)=g_0 (d\varphi (u), \, d\varphi (u))
\end{equation}
for each $u \in D_a$ (the tangent space to $D$ at $a$).  Since
$g_0$ is invariant under rotations around $0$,
\begin{equation}
  \label{eq:1.2}
  g_0 (z,z) = c|z|^2 \, , 
\end{equation}
where $c$ is a constant.  Here $D_0$ is identified with $\CC$.
Let $t \to z(t)$ be a curve with $z(0)=a$, $z'(0) =u \in \CC$.
Then $d\varphi (u)$ is the tangent vector
\begin{displaymath}
  \left\{ \frac{d}{dt} \varphi (z(t)) \right\}_{t=0}=
  \left( \frac{d \varphi}{dz}\right)_a
  \left( \frac{dz}{dt}\right)_0 =
\left\{ \frac{1-|a|^2}{(1-\bar{a}z)^2}\right\}_{z=a}\, u
\end{displaymath}
so by (\ref{eq:1.1}), (\ref{eq:1.2})
\begin{displaymath}
  g_a (u,u) =c \frac{1}{(1-|a|^2)^2} |u|^2 \, .
\end{displaymath}
Thus $g$ is the Riemannian structure
\begin{equation}
  \label{eq:1.3}
  ds^2 =c \frac{dx^2 + dy^2}{(1-x^2 -y^2)^2}
\end{equation}
and the proof shows that it is indeed invariant.

We shall now take $D$ as the unit disk $|z|<1$ with $g=ds^2$
given by (\ref{eq:1.3}) with $c=1$.  In our analysis on $D$ we are
mainly interested in the \emph{geodesics} in $D$ (the arcs
orthogonal to the boundary $B= \{ z \in \CC : |z| =1 \}$) and the
\emph{horocycles} in $D$ which are the circles inside $D$
tangential to $B$. Note that a horocycle tangential to $B$ at $b$
is orthogonal to all the geodesics in $D$ which end at $b$.

\section{The Non-Euclidean Fourier Transform}
\label{sec:2}

We first recall some of the principal results from Fourier
analysis on $\RR^n$.  The Fourier transform $f \to \tilde{f}$ or
$\RR^n$ is defined by
\begin{equation}
  \label{eq:2.1}
  \tilde{f} (u) = \int_{\RR^n} f(x) e^{-i(x,u)} \, dx
\end{equation}
where $(\, , \, )$ denotes the scalar product and $dx$ the
Lebesgue measure.  In polar coordinates $u=\lambda w$ $\lambda
\in \RR$, $w \in \SS^{n-1}$ we can write
\begin{equation}
  \label{eq:2.2}
  \tilde{f} (\lambda w) =\int_{\RR^n} f (x) e^{-i\lambda (x,w)} \,
  dx \, .
\end{equation}
It is then inverted by
\begin{equation}
  \label{eq:2.3}
  f(x) = (2\pi)^{-n} \int_{\RR^+ \times \SS^{n-1}}
      \tilde{f} (\lambda w)e^{i\lambda (x,w)}\lambda^{n-1}
      \, d\lambda \, dw
\end{equation}
say for $f \in \D (\RR^n)=\Cinf_c (\RR^n)$, $dw$ denoting the
surface element on the sphere $\SS^{n-1}$.  The Plancherel
formula
\begin{equation}
  \label{eq:2.4}
  \int_{\RR^n} |f(x)|^2 \, dx =(2\pi)^{-n}
           \int_{\RR^+ \times \SS^{n-1}} |\tilde{f} (\lambda ,w)|^2
           \lambda^{n-1} \, d\lambda \, dw
\end{equation}
expresses that $f \to \tilde{f}$ is an isometry of $L^2 (\RR^2)$
onto\break $L^2 (\RR^+ \times \SS^{n-1}, (2\pi)^{-n} \lambda^{n-1}\,
d\lambda \, dw)$.

The range of the mapping $f(x) \to \tilde{f} (\lambda w)$ as $f$
runs through $\D (\RR^n)$ is expressed in the following theorem
\cite{He7}.  A vector $a=(a_1 , \ldots ,a_n) \in \CC^n$ is said to
be \emph{isotropic} if $(a,a)=a^2_1 + \cdots + a^2_n=0$.

\begin{theorem}
  \label{th:2.1}

The Fourier transform $f(x) \to \tilde{f} (\lambda w)$ maps $\D
(\RR^n)$ onto the set of functions $\tilde{f} (\lambda w)
=\varphi (\lambda ,w) \in \Cinf (\RR \times \SS^{n-1})$
satisfying:

\romanparenlist

\begin{enumerate}
\item 
There exists a constant $A>0$ such that for each $w \in
\SS^{n-1}$ the function $\lambda \to \varphi (\lambda ,w)$
extends to a holomorphic function on $\CC$ with the property
\begin{equation}
  \label{eq:2.5}
\sup_{  \lambda \in \CC \, , \, w \in \SS^{n-1}}
    |\varphi (\lambda ,w)|(1+|\lambda|)^N 
    e^{-A |\Im \lambda }|< \infty
\end{equation}
for each $N \in \ZZ$.  ($\Im \lambda =$ imaginary part of
$\lambda$).

\item 
  For each $k \in \ZZ^+$ and each isotropic vector $a \in \CC^n$
  the function
  \begin{equation}
    \label{eq:2.6}
    \lambda \to \lambda^{-k} \int_{\SS^{n-1}}
    \varphi (\lambda ,w) (a,w)^k \, dw
  \end{equation}
is even and holomorphic in $\CC$.

\end{enumerate}

\end{theorem}

Condition~(\ref{eq:2.5}) expresses that the function $\lambda
\to \varphi (\lambda ,w)$ is of \emph{uniform exponential type}:
  The classical Paley--Wiener theorem states
  that $\D (\RR^n)^{\sim}$ consists of entire functions of
  exponential type in $n$ variables whereas in the description
  above only $\lambda$ enters.

Formula (\ref{eq:2.2}) motivates a Fourier transform definition on $D$.
The inner product $(x,\omega)$ equals the (signed) distance from $0$
to the hyperplane through $x$ with normal $\omega$.  A horocycle
in $D$ through $b$ is perpendicular to the (parallel) family
of geodesics ending at $b$ so is an analog of a hyperplane in
$\RR^n$.  Thus if $z \in D$, $b \in B$ we define $\langle z,b
\rangle$ as the (signed) distance from $0$ to the horocycle
through $z$ and $b$.  The Fourier transform $f \to \tilde{f}$ on
$D$ is thus defined by 
\begin{equation}
  \label{eq:2.7}
  \tilde{f} (\lambda ,b) = \int_{D} f(z)
  e^{(-i\lambda +1) \langle z,b \rangle}\, dz
\end{equation}
for all $b \in B$ and $\lambda \in \CC$ for which integral
converges.  Here $dz$ is the invariant surface element on $D$
\begin{equation}
  \label{eq:2.8}
  dz = (1-x^2 -y^2)^{-2} \, dx \, dy \, .
\end{equation}
The $+1$ term in \ref{eq:2.7} is included for later technical convenience.

The Fourier transform (\ref{eq:2.7}) is a special case of the
Fourier transform on a symmetric space $X=G/K$ of the non-compact
type, introduced in \cite{He3}.  Here $G$ is a semisimple connected
Lie group with finite center and $K$ is a maximal compact
subgroup.  In discussing the properties of $f \to \tilde{f}$
below we stick to the case $X=D$ for notational simplicity but
shall indicate (with references) the appropriate generalizations
to arbitrary $X$.  Some of the results require a rank restriction
on $X$.

\begin{theorem}
  \label{th:2.2}

The transform $f \to \tilde{f}$ in (\ref{eq:2.7}) is inverted by
\begin{equation}
  \label{eq:2.9}
  f(z) = \frac{1}{4\pi} \int_{\RR} \int_B
     \tilde{f} (\lambda ,b) e^{(i\lambda +1)\langle z,b\rangle}
     \lambda \thh\,  \left( \frac{\pi \lambda}{2}\right)
     \, d\lambda \, db \,, \quad f \in \D (D) \, .
\end{equation}
Also the map $f \to \tilde{f}$ extends to an isometry of $L^2
(D)$ onto $L^2(\RR^+ \times B,\mu)$ where $\mu$ is the measure 
\begin{equation}
  \label{eq:2.10}
    \mu =\frac{\lambda}{2\pi} \thh 
    \left( \frac{\pi \lambda}{2}\right)  \, d \lambda \, db
\end{equation}
and $db$ is normalized by $\int db =1$.
\end{theorem}

This result is valid for arbitrary $X=G/K$ (\cite{He3} \cite{He4}),
suitably formulated in terms of the fine structure of $G$.  While
this result resembles (\ref{eq:2.3})---(\ref{eq:2.4}) closely the
range theorem for $D$ takes a rather different form.

\begin{theorem}
  \label{th:2.3}
The range $\D (D)^{\sim}$ consists of the functions $\varphi
(\lambda ,b)$ which (in $\lambda$) are holomorphic of uniform
exponential type and satisfy the functional equation
\begin{equation}
  \label{eq:2.11}
  \int_B \varphi (\lambda ,b) e^{(i\lambda +1 \langle z,b  \rangle)}
  \, db = \int_B \varphi (-\lambda ,b)
  e^{(-i \lambda +1)\langle z,b \rangle} \, db \, .
\end{equation}

\end{theorem}

One can prove that condition (\ref{eq:2.11}) is equivalent to the
following conditions (\ref{eq:2.12}) for the Fourier coefficients
$\varphi_k (\lambda)$ of $\varphi$
\begin{eqnarray}
\nonumber
  \varphi_k (\lambda) &=& \frac{1}{2\pi}\int^{2\pi}_0
     \varphi (\lambda ,e^{i\theta}) e^{-ik\theta}\, d\theta\\
\label{eq:2.12}
  \varphi_k (-\lambda) p_k (-i\lambda) &=&
     \varphi_k (\lambda)p_k (i\lambda )\, ,\quad k \in \ZZ \, ,
\end{eqnarray}
where $p_k (x)$ is the polynomial
\begin{displaymath}
  p_k (x) =\frac{\Gamma (\frac{1}{2}(x+1)+|k|)}
       {\Gamma (\frac{1}{2}(x+1))}\, .
\end{displaymath}
Again these results are valid for arbitrary $X=G/K$ (\cite{He5}
and \cite{He7}).

The Paley--Wiener type theorems can be extended to the Schwartz
spaces $\s^p(D)$ $(0<p\leq 2)$.  Roughly speaking, $f $ belongs
to $\s^p(D)$ if each invariant derivative $Df$ belongs to
$L^p(D)$, more precisely, it is rapidly decreasing in the
distance from $0$ even after multiplication by the $p$\st{th}
root of the volume element.  Let $S_p$ denote the strip $|\Im
\lambda | <\frac{2}{p}-1$ in $\CC$ and $\s (S_p \times B)$ the
space of smooth functions on $\s (S_p \times B)$ holomorphic (in
$\lambda$) in $S_p$ and rapidly decreasing (uniformly for $b \in
B$) on each line $\lambda =\xi + i\eta$ $(|\eta
|<\frac{2}{p}-1)$.  Then we have

\begin{theorem}
  \label{th:2.4}

The Fourier transform $f \to \tilde{f}$ on $D$ is a bijection of
$\s^p (D)$ onto the set of $\varphi \in \s (S_p \times B)$
satisfying (\ref{eq:2.11}).
\end{theorem}

The theorem holds for all $X=G/K$ (Eguchi \cite{Eg}).  The proof is
complicated.  For the case of $K$-invariant functions (done for
$p=2$ by Harish--Chandra \cite{H} and Trombi--Varadarajan \cite{TV}
for general $p$) a substantial simplification was done by Anker
\cite{A}.  A further range theorem for the space of functions for
which each invariant derivative has arbitrary exponential decay
was proved by Oshima, Saburi and Wakayama \cite{OSW}.  See also
Barker \cite{Bar} (p.~27) for the operator Fourier transform of the
intersection of all the Schwartz spaces on $\SL (2,\RR)$.

In classical Fourier analysis on $\RR^n$ the Riemann--Lebesgue
lemma states that for $f \in L^1 (\RR)$,   $\tilde{f}$ tends to $0$ at $\infty$.  For $D$
the situation is a bit different.

\begin{theorem}
  \label{th:2.5}
Let $f \in L^1(D)$.  Then there exists a null set $N$ in $B$ such
that if $b \in B-N$, $\lambda \to\tilde{f} (\lambda ,b)$ is
holomorphic in the strip $|\Im \lambda|<1$ and
\begin{equation}
  \label{eq:2.13}
  \lim_{\xi \to \infty} \tilde{f} (\xi +i\eta ,b)=0
\end{equation}
uniformly for $|\eta|\leq 1$.

\end{theorem}

The proof \cite{HRSS} is valid even for symmetric spaces $X=G/K$ of
arbitrary rank.  Moreover
\begin{equation}
  \label{eq:2.14}
  \| \tilde{f} (\lambda ,\cdot)\|_1 \to 0
  \hbox{ as } \lambda \to \infty\, , 
\end{equation}
uniformly in the strip $|\Im \lambda| \leq 1$, and this extends
to $f \in L^p\,\, (1 \leq p<2)$ this time in the strip $|\Im
\lambda| <\frac{2}{p}-1$ (\cite{SS}).  In particular, if $f \in L^p
(D)$ then there is a null set $N$ in $B$ such that $\tilde{f}
(\lambda ,b)$ exists for $b\notin N$ and all $\lambda$ in the
strip $|\Im \lambda| <\frac{2}{p}-1$.

The classical inversion formula for the Fourier transform on
$\RR^n$ now extends to $f \in L^p (D)$ $(1 \leq p<2)$ as follows.

\begin{theorem}
  \label{th:2.6}
Let $f \in L^p(D)$ and assume $\tilde{f} \in L^1 (\RR \times B
,\mu)$ (with $\mu$ as in (\ref{eq:2.10})).  Then the inversion
formula (\ref{eq:2.9}) holds for almost all $z \in D$ (the
Lebesgue set for $f$).

\end{theorem}

Again this holds for all $X=G/K$.  A result of this type was
proved by Stanton and Thomas \cite{ST} without invoking $\tilde{f}$
explicitly (since the existence had not been established).  The
version in Theorem~\ref{th:2.6} is from \cite{SS}.

In Schwartz's theory of mean--periodic functions \cite{Sc} it is
proved that any closed translation--invariant subspace of $\Cinf
(\RR)$ contains an exponential $e^{\mu x}$.  The analogous
question here would be:

\emph{Does an arbitrary closed invariant subspace of $\Cinf (D)$
  contain an exponential}
\begin{equation}
  \label{eq:2.15}
  e_{\mu ,b} (z)=  e^{\mu \langle z,b \rangle}
\end{equation}
\emph{for some }$\mu \in \CC$ \emph{and some }$b \in B$\emph{?}

Here the topology of $\Cinf (D)$ is the usual Fr\'echet space
topology and ``invariant'' refers to the action of the group
$G=\SU (1,1)$ on $D$.  The answer is yes.

\begin{theorem}
  \label{th:2.7}

Each closed invariant subspace $E$ of $\Cinf (D)$ contains an
exponential $e_{\mu ,b}$.
\end{theorem}

This was proved in \cite{HS} for all symmetric $X=G/K$ of rank
one.  Here is a sketch of the proof.  By a result of Bagchi and
Sitaram \cite{BS} $E$ contains a spherical function
\begin{equation}
  \label{eq:2.16}
  \varphi_{\lambda} (z) =\int_B
     e^{(i\lambda +1)\langle z,b \rangle} \, db \, ,
     \quad \varphi_{\lambda} =\varphi_{-\lambda} \, .
\end{equation}
For either $\lambda$ or $-\lambda$ it is true (\cite{He9},
Lemma~2.3, Ch.~III) that the Poisson transform $P_{\lambda}:F \to
f$ where 
\begin{equation}
  \label{eq:2.17}
  f(z) =\int_B e^{(i\lambda +1)\langle z,b \rangle}
  F(b) \, db \, , 
\end{equation}
maps $L^2 (B)$ into the closed invariant subspace of $E$
generated by $\varphi_{\lambda}$.  On the other hand it is
proved in \cite{He9} (Ex.~B1 in Ch.~III) that $e_{i\lambda +1,b}$
is a  series of terms $P_{\lambda} (F_n)$ where $F_n
\in L^2 (B)$ and the series converges in the topology of $\Cinf (D)$.  Thus $e_{i\lambda +1,b} \in E$ as desired.

The following result for the Fourier transform on $\RR^n$ is
closely related to the Wiener Tauberian theorem.

\emph{Let} $f \in L^1(\RR^n)$ \emph{ be such that }$\tilde{f} (u) \neq
0$ \emph{for all }$u \in \RR^n$.  \emph{Then the translates of
}$f$ \emph{span a dense subspace  of }$L^1(\RR^n)$.

There has been considerable activity in establishing analogs of
this theorem for semisimple Lie groups and symmetric spaces.  See
e.g.~\cite{EM}, \cite{Sa}, \cite{Si1}, \cite{Si2}.  The neatest version for
$D$ seems to me to be the following result from \cite{SS} \cite{MRSS}
which remains valid for $X=G/K$ of rank one.

Let $d (z,w)$ denote the distance in $D$ and if $\epsilon >0$, let
$\L_{\epsilon}(D)$ denote the space of measurable functions $f$
on $D$ such that $\int_D |f(z)| e^{\epsilon d (0,z)} \, dz
<\infty$.  Let $T_{\epsilon}$ denote the strip $|\Im \lambda|\leq
1+\epsilon$.

\begin{theorem}
  \label{th:2.8}
Let $f \in \L_{\epsilon}(X)$ and assume $f$ is not almost
everywhere equal to any real analytic function.  Let 
\begin{displaymath}
Z=\{ \lambda \in T_{\epsilon} : \tilde{f} (\lambda ,\cdot)\equiv 0\}\,.
\end{displaymath}
If $Z =\emptyset$ then the translates of $f$ span a dense
subspace of $L^1(D)$.
\end{theorem}

A theorem of Hardy's on Fourier transforms on $\RR^n$ asserts in a
precise fashion that $f$ and its Fourier transform cannot both
vanish too fast at infinity.  More precisely (\cite{Ha}):

Assume
\begin{displaymath}
  |f(x)|\leq A e^{-\alpha |x|^2}\, , \, 
  |\tilde{f} (u)|\leq B e^{-\beta |u|^2}\, , 
\end{displaymath}
where $A$, $B$, $\alpha$ and $\beta$ are positive constants and
$\alpha \beta >\tfrac{1}{4}$.  Then $f \equiv 0$.

Variations of this theorem for $L^p$ spaces have been proved by
Morgan \cite{M} and Cowling--Price \cite{CP}.

For the Fourier transform on $D$ the following result holds.

\begin{theorem}
  \label{th:2.9}

Let $f$ be a measurable function on $D$ satisfying
\begin{displaymath}
  |f(x)| \leq  C e^{-\alpha d (0,x)^2}\quad
  |\tilde{f}(\lambda ,b)| \leq  C e^{-\beta |\lambda|^2}
\end{displaymath}
where $C$, $\alpha$, $\beta$ are positive constants.  If $\alpha
\beta >16$ then $f \equiv 0$.
\end{theorem}

This is contained in Sitaram and Sundari \cite{SiSu} \S~5 where an
extension to certain symmetric spaces $X=G/K$ is also proved.
The theorem for all such $X$ was obtained by Sengupta \cite{Se},
together with refinements in terms of $L^p(X)$.

Many such completions of Hardy's theorem have been given, see
\cite{RS}, \cite{CSS}, \cite{NR}, \cite{Shi}.

\section{Eigenfunctions of the Laplacian}
\label{sec:3}

Consider first the plane $\RR^2$ and the Laplacian 
\begin{displaymath}
  L^0 =\frac{\partial^2}{\partial x^2_1} +
       \frac{\partial^2}{\partial x^2_2} \, .
\end{displaymath}
Given a unit vector $\omega \in \RR^2$ and $\lambda \in \CC$ the
function $x \to e^{i\lambda (x,\omega)}$ is an eigenfunction
\begin{equation}
  \label{eq:3.1}
  L^0_x e^{i\lambda (x,\omega)}=-\lambda^2 e^{i\lambda (x,\omega)} \, .
\end{equation}
Because of (\ref{eq:2.3}) one might expect all eigenfunctions of
$L$ to be a ``decomposition'' into such eigenfunctions with fixed
$\lambda$ but variable $\omega$.

Note that the function $\omega \to e^{i\lambda (x,\omega)}$ is
the restriction to $\SS^1$ of the holomorphic function
\begin{displaymath}
  z \to \exp \left[ \tfrac{1}{2} (i\lambda )x_1
       (z+z^{-1})+\tfrac{1}{2}\lambda x_2 (z-z^{-1})\right] \, 
     \quad z \in \CC -(0) \, ,
\end{displaymath}
which satisfies a condition
\begin{equation}
  \label{eq:3.2}
  \sup_z \left( |f(z)| e^{-a|z|-b|z|^{-1}}\right)
     <\infty \, , 
\end{equation}
with $a,b \geq 0$.  Let $E_{a,b}$ denote the Banach space of
holomorphic functions satisfying (\ref{eq:3.2}), the norm being
the expression in (\ref{eq:3.2}).  We let $E$ denote the union of
the spaces $E_{a,b}$ and give it the induced topology.  We
identify the elements of $E$ with their restrictions to $\SS^1$
and call the members of the dual space $E'$ \emph{entire functionals}.

\begin{theorem}(\cite{He6})
  \label{th:3.1}
The eigenfunctions of $L^0$ on $\RR^2$ are precisely the harmonic
functions and the functions
\begin{equation}
  \label{eq:3.3}
  f(x) =\int_{\SS^1} e^{i\lambda (x,\omega)} \, dT (\omega)
\end{equation}
where $\lambda \in \CC -(0)$ and $T$ is an entire functional on $\SS^1$.
\end{theorem}

For the non-Euclidean metric (\ref{eq:1.3}) (with $c=1$) the
Laplacian is given by
\begin{equation}
  \label{eq:3.4}
  L=(1-x^2-y^2)^2 \left( \frac{\partial^2}{\partial x^2} 
       + \frac{\partial^2}{\partial y^2} \right)
\end{equation}
and the exponential function $e_{\mu,b} (z)=e^{\mu \langle z,b
    \rangle}$ is an eigenfunction:
\begin{equation}
  \label{eq:3.5}
  L_z (e^{(i\lambda +1)\langle z,b \rangle})
     =-(\lambda^2 +1) e^{(i\lambda +1)\langle z,b \rangle}\, .
\end{equation}
In particular, the function $z \to e^{2\langle z,b \rangle}$ is a
harmonic function and in fact coincides with the classical
Poisson kernel from potential theory:
\begin{equation}
  \label{eq:3.6}
  e^{2\langle z,b \rangle}=\frac{1-|z|^2}{|z-b|^2}\, .
\end{equation}
Again the eigenfunctions of $L$ are obtained from the functions
$e_{\mu ,b}$ by superposition.  To describe this precisely
consider the space $\A (B)$ of analytic functions on $B$.  Each
$F \in \A (B)$ extends to a holomorphic function on a belt
$B_{\epsilon} :1-\epsilon <|z|<1+\epsilon$ around $B$.  The
space $\Hs (B_{\epsilon})$ of holomorphic functions on
$B_{\epsilon}$ is topologized by uniform convergence on compact
subsets.  We can view $\A (B)$ as the union $\cup^{\infty}_n \Hs
(B_{1/n})$ and give it the inductive limit topology.  
The dual space $\A'(B)$ then consists of the \emph{analytic
  functionals} on $B$ (or hyperfunctions in $B$).

\begin{theorem}(\cite{He4}, IV, \S1).  
  \label{th:3.2}
The eigenfunctions of $L$ are precisely the functions
\begin{equation}
  \label{eq:3.7}
  u(z) =\int_B e^{\mu \langle z ,b \rangle}\, dT (b) \, , 
\end{equation}
where $\mu \in \CC$ and $T \in \A' (B)$.
\end{theorem}

Lewis in \cite{L} has proved (under minor restriction on $\mu$)
that $T$ in (\ref{eq:3.7}) is a distribution if and only if $u$
has at most an exponential growth (in $d(0,z)$).  On the other
hand, Ban and Schlichtkrull proved in \cite{BaS} that $T \in \Cinf
(B)$ if and only if all the invariant derivatives of $u$ have the
same exponential growth.

We consider now the natural group representations on the
eigenspaces.  The group $\MM (2)$ of isometries of $\RR^2$ acts
transitively on $\RR^2$ and leaves the Laplacian $L^0$ invariant:
 $L^0 (f \circ \tau)=(L^0f)\circ \tau$ for each $\tau \in \MM
 (2)$.  If $\lambda \in \CC$ the eigenspace 
 \begin{displaymath}
   \E_{\lambda} = \{ f \in \Cinf (\RR^2) : 
   L^0 f=-\lambda^2 f \}
 \end{displaymath}
is invariant under the action $f \to f\circ \tau^{-1}$ so we have
a representation $T_{\lambda}$ of $\MM (2)$ on $\E_{\lambda}$
given by $T_{\lambda} (\tau)(f) = f\circ \tau^{-1}$, the
\emph{eigenspace representation}.

\begin{theorem}({\cite{He6}})
  \label{th:3.3}
The representation $T_{\lambda}$ is irreducible if and only if
$\lambda \neq 0$.
\end{theorem}

Similarly the group $G=\SU (1,1)$ of conformal transformations
\begin{displaymath}
  z \to \frac{az +b}{\bar{b} z+a} \quad (|a|^2 -|b|^2 =1)
\end{displaymath}
leaves (\ref{eq:1.3}) and the operator $L$ in (\ref{eq:3.4})
invariant.  Thus we get again an eigenspace representation
$\tau_{\lambda}$ of $G$ on each eigenspace
\begin{displaymath}
  \E_{\lambda} = \{ f \in \Cinf (D) : Lf =-(\lambda^2 +1)f\} \, .
\end{displaymath}

\begin{theorem}({\cite{He4}})
  \label{th:3.4}
The representation $\tau_{\lambda}$ is irreducible if and only if
$i\lambda +1 \notin 2 \ZZ$.
\end{theorem}

Again all these results extend to Euclidean spaces of higher
dimensions and suitably formulated, to all symmetric spaces $G/K$
of the noncompact type.

\section{The Radon Transform}
\label{sec:4}

\subsection*{\large{A.  The Euclidean Case.}}

Let $d$ be a fixed integer, $0<d<n$ and let $\GG (d,n)$ denote
the space of $d$-dimensional planes in $\RR^n$.  To a function
$f$ on $\RR^n$ we associate a function $\hat{f}$ on $\GG (d,n)$
by
\begin{equation}
  \label{eq:4.1}
  \hat{f} (\xi) =\int_{\xi} f(x) \, dm (x) \, , \quad
  \xi \in \GG (d,n) \, ,
\end{equation}
$dm$ being the Euclidean measure on $\xi$.  The transform $f \to
\hat{f}$ is called the \emph{$d$-plane transform}. For $d=1$,
$n=2$ it is the classical Radon transform. The parity of
$d$ turns out to be important.

The inversion of the transform $f \to \hat{f}$ is well known
(case $d=n-1$ in \cite{R}, \cite{J}, \cite{GS}, general $d$ in \cite{F}, \cite{He1},
\cite{He2}).  We shall give another group-theoretic method here,
resulting in alternative inversion formulas.

The group $G=\MM (n)$ acts transitively both on $\RR^n$ and on
$\GG (d,n)$.  In particular, $\RR^n =G/K$ where $K=\OO (n)$.  Let
$p>0$.  Consider a pair $x \in \RR^n$, $\xi \in \GG (d,n)$ at
distance $p=d(x,\xi)$.  Let $g \in G$ be such that $g \cdot
0=x$.  Then the family $kg^{-1}\cdot \xi$ constitutes the set of
elements in $\GG (d,n)$ at distance $p$ from $0$.  Along with the
transform $f \to \hat{f}$ we consider the \emph{dual transform}
$\varphi \to \sucheck{\varphi}$ given by
\begin{equation}
  \label{eq:4.2}
  \sucheck{\varphi} (x) = \int_{\xi \ni x} \varphi (\xi) \, 
  d \mu (\xi) \, , 
\end{equation}
the average of $\varphi$ over the set of $d$-planes passing
through $x$.  More generally we put
\begin{equation}
  \label{eq:4.3}
  \sucheck{\varphi}_p (x) =\int_{d (\xi ,x)=p} 
     \varphi (\xi) \, d\mu (\xi)
\end{equation}
the average of $\varphi$ over the set of $d$-planes at distance
$p$ from $x$.  Since $K$ acts transitively on the set of
$d$-planes through $0$ we see by the above that
\begin{equation}
  \sucheck{\varphi}_p (g \cdot 0) =\int_K \varphi 
     (gkg^{-1} \cdot \xi) \, dk \, ,
\end{equation}
$dk$ being the normalized Haar measure on $K$.  Let $(M^rf)(x)$
denote the mean-value of $f$ over the sphere $S_r(x)$ of radius
$r$ with center $x$.  If $z \in
\RR^n$ has distance $r$ from $0$ we then have
\begin{equation}
  \label{eq:4.5}
  (M^rf)(g \cdot 0)=\int_K f(gk\cdot z)\, dk \, .
\end{equation}
We thus see that since $d(0,g^{-1}\cdot y)=d (x,y)$,
\begin{eqnarray*}
  (\hat{f})\spcheck_p (x) &=& \int_K \hat{f}(gkg^{-1}\cdot \xi)\, dk
     = \int_{\xi}\, dk \int_{\xi} f(gkg^{-1}\cdot y)\, dm (y)\\
     &=& \int_{\xi} \, dm (y) \int_K f(gkg^{-1} \cdot y)\, dk
     = \int_{\xi} (M^{d(x,y)}f)(x) \, dm (y) \, .
\end{eqnarray*}
Let $x_0$ be the point in $\xi$ at minimum distance $p$ from
$x$.  The integrand $(M^{d(x,y)}f(x))$ is constant in $y$ on each
sphere in $\xi$ with center $x_0$.  It follows that
\begin{equation}
  \label{eq:4.6}
  (\hat{f})\spcheck_p (x) = \Omega_d \int^{\infty}_0
   (M^qf) (x) r^{d-1}\, dr \, , 
\end{equation}
where $r=d(x_0,y)$, $q=d(x,y)$, $\Omega_d$ denoting the area of
the unit sphere in $\RR^d$.  We have $q^2=p^2+r^2$ so putting
$F(q) =(M^qf) (x)$, $\hat{F} (p) =(\hat{f})\spcheck_p (x)$ we
have
\begin{equation}
  \label{eq:4.7}
  \hat{F} (p) = \Omega_d \int^{\infty}_p
     F(q) (q^2-p^2)^{d/2-1} q \, dq \, .
\end{equation}
This Abel-type integral equation has an inversion
\begin{equation}
  \label{eq:4.8}
  F(r) =c(d) \left( \frac{d}{d(r^2)}\right)^d
    \int^{\infty}_r p(p^2-r^2)^{d/2-1}\hat{F} (p) \, dp \, , 
\end{equation}
where $c(d)$ is a constant, depending only on $d$.  Putting $r=0$
we obtain the inversion formula
\begin{equation}
  \label{eq:4.9}
  f(x) = c(d) \left[ \left( \frac{d}{d(r^2)}\right)^d
       \int^{\infty}_r p (p^2-r^2)^{d/2-1}(\hat{f})\spcheck_p
         (x) \, dp \right]_{r=0} \, .
\end{equation}

Note that in (\ref{eq:4.8})
\begin{displaymath}
  p(p^2-r^2)^{d/2-1} = \frac{d}{dp}(p^2-r^2)^{(d/2)}\cdot
     \frac{1}{d}
\end{displaymath}
so in (\ref{eq:4.8}) we can use integration by parts and the
integral becomes
\begin{displaymath}
  C \int^{\infty}_r (p^2-r^2)^{d/2} \hat{F}' (p) \, dp \, .
\end{displaymath}
Applying $\displaystyle{\frac{d}{d(r^2)} = \frac{1}{2r} \,
  \frac{d}{dr}}$ to this integral reduces the exponent $d/2$ by
$1$.  For $d$ odd we continue the differentiation
$\tfrac{d+1}{2}$ times until the exponent is $-\tfrac{1}{2}$.
For $d$ even we continue until the exponent is $0$ and then
replace $\int^{\infty}_r \hat{F}' (p) \, dp$ by $-\hat{F}(r)$.
 This $\hat{F} (r)$ is an even function so taking
 $(d/d(r^2))^{d/2}$ at $r=0$ amounts to taking a constant
 multiple of $(d/dr)^d$ at $r=0$.  We thus get the following
 refinement of (\ref{eq:4.9}) where we recall that
 $(\hat{f})\spcheck_r(x)$ is the average of the integrals of $f$
 over the $d$-planes tangent to $S_r (x)$.

 \begin{theorem}
\label{th:4.1}
The $d$-plane transform is inverted as follows:

\romanparenlist
\begin{enumerate}
\item 
If $d$ is even then
\begin{equation}
  \label{eq:4.10}
  f(x) = C_1 \left[ \left( \frac{d}{dr}\right)^d (\hat{f})\spcheck_r
         (x) \right]_{r=0}\, .
\end{equation}

\item 
If $d$ is odd then
\begin{equation}
   \label{eq:4.11}
 \hspace*{-7ex}  f(x) = C_2 \left[ \left( \frac{d}{d(r^2)}\right)^{(d-1)/2}
          \int^{\infty}_r (p^2-r^2)^{-1/2}\frac{d}{dp}
          (\hat{f})\spcheck_p (x) \, dp \right]_{r=0}\, . 
\end{equation}

\item 
  If $d=1$ then
  \begin{equation}
    \label{eq:4.12}
    f(x) =-\frac{1}\pi{\int^{\infty}_0}\frac{1}{p}\, 
    \frac{d}{dp}(\hat{f})\spcheck_p (x) \, dp \, .
  \end{equation}

\end{enumerate}
 \end{theorem}

For $n=2$ formula (\ref{eq:4.12}) is proved in Radon's original
paper \cite{R}.  Note that the constant $-1/\pi$ is the same for
all $n$.  In the case $d=n-1$ the formula in (i) coincides with
formula (21) in Rouvi\`ere \cite{Ro}.

Another inversion formula (\cite{He1}, \cite{He2}) valid for all
$d$ and $n$ is
\begin{eqnarray}
  \label{eq:4.13}
  f &=& c (-L)^{d/2} ((\hat{f})\spcheck) \, \\
\noalign{\nonumber\hbox{where}}\\
\nonumber
  c &=& \frac{\Gamma \left( \tfrac{n-d}{2}\right)}
       {(4\pi)^{d/2}\Gamma \left( \frac{n}{2}\right)}\, .
\end{eqnarray}
Here the fractional power of $L$ is defined in the usual way by
the Fourier transform.  The parity of $d$ shows up in the same
way as in Theorem~\ref{th:4.1}.

For range questions for the transform $f \to \hat{f}$ see an
account in \cite{He10} and references there.

\subsection*{\large{B.  The Hyperbolic Case.}}

The hyperbolic space $\HH^n$ is the higher-dimensional version of
(\ref{eq:1.3}) and its Riemannian structure is given by
\begin{equation}
  \label{eq:4.14}
  ds^2 =4 \frac{dx^2_1 + \cdots + dx^2_n}{(1-x^2_1 -\cdots - x^2_n)^2}
\end{equation}
in the unit ball $|x|<1$.  The constant~$4$ is chosen such that
the curvature is now $-1$.  The $d$-dimensional totally geodesic
submanifolds are spherical caps perpendicular to the boundary $B:
|x|=1$.  They are natural analogs of the $d$-planes in $\RR^n$.
We have accordingly a Radon transform $f \to \hat{f}$, where
\begin{equation}
  \label{eq:4.15}
  \hat{f} (\xi) =\int_{\xi} f(x) \, dm (x) \, \quad
  \xi \in \Xi \, ,
\end{equation}
where $\Xi$ is the space of $d$-dimensional totally geodesic
submanifolds of $\HH^n$.

The group $G$ of isometries of $\HH^n$ acts transitively on
$\Xi$ as well.  As in (\ref{eq:4.2})---(\ref{eq:4.3}) we
consider the dual transform
\begin{equation}
\label{eq:4.16}
\sucheck{\varphi} (x) = \int_{\xi \ni x} \varphi (\xi)\, 
   d \mu (\xi)
\end{equation}
and more generally for $p \geq 0$\, , 
\begin{equation}
  \label{eq:4.17}
  \sucheck{\varphi}_p (x) = \int_{d(\xi ,x)=p}
     \varphi (\xi) \, d\mu (\xi) \, , 
\end{equation}
the mean value of $\varphi $ over the set of $\xi \in \Xi$ at
distance $p$ from $x$.  The formula
\begin{equation}
  \label{eq:4.18}
  (\hat{f})\spcheck_p (x) = \int_{\xi} (M^{d(x,y)}f)
     (x) \, dm (y)
\end{equation}
is proved just as before.  Let $x_0$ be the point in $\xi$ at
minimum distance $p$ from $x$ and put $r =d (x_0 ,y)$,
$q=d(x,y)$.  Since the geodesic triangle $(xx_0 y)$ is right
angled at $x_0$ we have by the cosine rule
\begin{equation}
  \label{eq:4.19}
  \cosh q = \cosh p \, \cosh r \, .
\end{equation}
Also note that since $\xi$ is totally geodesic, distances between
two points in $\xi$ are the same as in $\HH^n$.  In particular
$(M^{d(x,y)}f)(x)$ is constant as $y$ varies on a sphere in $\xi$
with center $x_0$.  Therefore  (\ref{eq:4.18}) implies
\begin{equation}
  \label{eq:4.20}
  (\hat{f})\spcheck_p (x) =\Omega_d \int^{\infty}_0
     (M^q f) (x) \sinh^{d-1} r \, dr \, .
\end{equation}
For $x$ fixed we put
\begin{displaymath}
F (\cosh q) = (M^q f) (x) \, , \quad
\hat{F} (\cosh p) =(\hat{f})\spcheck_p (x) \, , 
\end{displaymath}
substitute in (\ref{eq:4.20}) and use (\ref{eq:4.19}).  Writing
$t=\cosh p$,  $s=\cosh r$ we obtain the integral equation
\begin{equation}
  \label{eq:4.21}
  \hat{F} (t) = \Omega_d \int^{\infty}_1 F (ts)(s^2-1)^{d/2 \,
    -1} \, ds \, .
\end{equation}
Putting here $u=ts, \quad ds = t^{-1}\, du$ we get the Abel--type
integral equation
\begin{displaymath}
  t^{d-1} \hat{F} (t) = \Omega_d \int^{\infty}_t
     u^{-1} F(u) (u^2-t^2)^{d/2 \, -1} u \, du \, ,
\end{displaymath}
which by (\ref{eq:4.8}) is inverted by
\begin{equation}
  \label{eq:4.22}
  r^{-1} F(r) =c(d) \left( \frac{d}{d(r^2)}\right)^d
     \int^{\infty}_r t(t^2-r^2)^{d/2 \, -1}
     t^{d-1} \hat{F} (t) \, dt \, .
\end{equation}
Here we put $r=1$ and $s(p) =\cosh^{-1} p$.  We then obtain the
following variation of Theorem~3.12, Ch.~I in \cite{He9}:

\begin{theorem}
  \label{th:4.2}
The transform $f \to \tilde{f}$ is inverted by
\begin{equation}
  \label{eq:4.23}
  f(x) = C \left[\left( \frac{d}{d(r^2)} \right)^d
     \int^{\infty}_r (t^2-r^2)^{d/2 \, -1}t^d
     (\hat{f})\spcheck_{s(t)} (x) \, dt \right]_{r=1 \, .}
\end{equation}

\end{theorem}

As in the proof of Theorem~\ref{th:4.1} we can obtain the
following improvement.

\begin{theorem}
  \label{th:4.3}
\romanparenlist
\begin{enumerate}
\item 
If $d$ is even the inversion can be written
\begin{equation}
  \label{eq:4.24}
  f(x) =C \left[ \left( \frac{d}{d(r^2)}\right)^{d/2}
       (r^{d-1}(\hat{f})\spcheck_{s(r)}(x) ) \right]_{r=1}\, .
\end{equation}

\item 
If $d=1$ then
\begin{equation}
  \label{eq:4.25}
  f(x) =-\tfrac{1}{\pi} \int^{\infty}_0 \, \frac{1}{\sinh p}
     \frac{d}{dp} \left( (\hat{f})\spcheck_p (x)\right)\, dp \, .
\end{equation}

\end{enumerate}
\end{theorem}

\begin{proof}
  Part (i) is proved as (\ref{eq:4.10}) except that we no longer
  can equate $(d/d (r^2))^{d/2}$ with $(d/dr)^d$ at $r=1$.

For (ii) we deduce from (\ref{eq:4.22}) since $t(t^2-r^2)^{-1/2}
=\tfrac{d}{dt} (t^2-r^2)^{1/2}$ that
\begin{eqnarray*}
  F(1) &=& - \frac{c(1)}{2} \left[ \frac{d}{dr}\int^{\infty}_r
    (t^2-r^2)^{1/2} \frac{d}{dt}\hat{F} (t) \, dt \right]_{r=1}\\
  &=& \frac{c(1)}{2}\int^{\infty}_1 (t^2 -1)^{-1/2}
      \frac{d}{dt}\hat{F} (t) \, dt \, .
\end{eqnarray*}
Putting again $t=\cosh p$, $dt=\sinh p \, dp$ our expression
becomes
\begin{displaymath}
  \frac{c(1)}{2} \int^{\infty}_0 \frac{1}{\sinh p} \,
    \frac{d}{dp} ((\hat{f})\spcheck_p)(x) \, dp \, .
\end{displaymath}
\end{proof}

\begin{remark}
  For $n=2$, $d=1$ formula (\ref{eq:4.25}) is stated in Radon
    \cite{R}, Part~C.  The proof (which is only indicated) is
    very elegant but would not work for $n>2$.

For $d$ even (\ref{eq:4.24}) can be written in a simpler form
(\cite{He1}) namely
\begin{equation}
  \label{eq:4.26}
  f=c \,\, Q_d (L) ((\hat{f})\spcheck)\, ,
\end{equation}
where $\displaystyle{c=\tfrac{\Gamma \left(
    \tfrac{n-d}{2}\right)}{(-4\pi)^{d/2}\Gamma \left(
      \tfrac{n}{2}\right)}}$ and $Q_d$ is the polynomial
\begin{displaymath}
  Q_d (x) = (x+(d-1)(n-d)) (x+(d-3)(n-d+2))\cdots
     (x+1 \cdot (n-2))\, .
\end{displaymath}

The case $d=1$, $n=2$ is that of the X-ray transform on the
non-Euclidean disk ((\ref{eq:4.15}) for $n=2$).  Here are two
further alternatives to the inversion formula (\ref{eq:4.25}).
Let $S$ denote the integral operator
\begin{equation}
  \label{eq:4.27}
  (Sf)(x) = \int_D (\coth \, d(x,y)-1) f(y) \, dy\, .
\end{equation}
Then
\begin{equation}
  \label{eq:4.28}
  LS (\hat{f})\spcheck =-4\pi^2 f \, , \quad
     f \in \D (X) \, .
\end{equation}
This is proved by Berenstein--Casadio \cite{BC}; see \cite{He10}
for a minor simplification.  By invariance it suffices to prove
(\ref{eq:4.28}) for $f$ radial and then it is verified by taking
the spherical transform on both sides.  Less explicit versions of
(\ref{eq:4.28}) are obtained in \cite{BC} for any dimension $n$
and $d$.

\end{remark}

One more inversion formula was obtained by Lissianoi and
Ponomarev \cite{LP} using (\ref{eq:4.23}) for $d=1$, $n=2$ as a
starting point.  By parameterizing the geodesics $\gamma$ by the
two points of intersection of $\gamma$ with $B$ they prove a
hyperbolic analog of the Euclidean formula:
\begin{equation}
 \label{eq:30}
f(x) = \int_{\SS^1} \left\{ \Hs_p \frac{d}{dp}\hat{f} (\omega,p)
     \right\}_{p=(\omega ,x)}\, d\omega \, ,
\end{equation}
which is an alternative to (\ref{eq:4.12}).  Here $\Hs_p$ is a
normalized Hilbert transform in the variable $p$ and $\hat{f}
(\omega ,p)$ is $\hat{f} (\xi)$ for the line $(x,\omega)=p $, where
$|\omega |=1$.

In the theorems in this section we have not discussed smoothness
and decay at infinity of the functions.  Here we refer to
\cite{Je}, \cite{Ru1},  \cite{Ru2},  \cite{BeR1} and
\cite{BeR2}as examples.

Additional inversion formulas for the transform $f \to \hat{f}$
can be found in \cite{Sem}, \cite{Ru3} and \cite{K}.  The range
problem for the transform $f \to \hat{f}$ is treated in
\cite{BCK} and \cite{I}.

\subsection*{Added in Proof:}

I have since this was written proved an inversion formula for the
X-ray transform on a noncompact symmetric space of rank $l > 1$.  It is
similar to~(4.12) except that in~(4.17) one restricts the averaging  to the
set of geodesics each of which lies in a flat $l$-dimensional totally
geodesic submanifold through $x$ and at distance $p$ from $x$.  On the other
hand,   Rouvi\`ere had proved earlier that the inversion formula~(4.25) holds
almost unchanged for the X-ray transform on a noncompact symmetric space
of rank~$1$.

\end{document}